\documentclass[12pt]{article}

\usepackage[T2A]{fontenc}
\usepackage[utf8]{inputenc}

\usepackage{amsfonts}
\usepackage{amssymb}
\usepackage{amsmath}
\usepackage{amsthm}

\def \le {\leqslant}
\def \ge {\geqslant}

\begin{document}

\begin{Large}
 \centerline{\bf On Minkowski diagonal continued fraction}
\end{Large}
\vskip+0.5cm

\centerline{ by {\bf Nikolay Moshchevitin}\footnote{research is supported by RFBR grant No.12-01-00681-a}}

\vskip+0.5cm
\begin{small}
 {\bf Abstract.}\, We study some properties of the function $\mu_\alpha (t)$ associated with the 
Minkowski diagonal continued fraction for real $\alpha$.
\end{small}
\vskip+0.5cm

{\bf 1. Irrationality measure function, Lagrange and Dirichlet spectra.}

For a real $\alpha$ we consider  the irrationality measure function
$$
\psi_\alpha (t) = \min_{1\le x \le t,\, x \in \mathbb{Z}}||x\alpha ||
$$ 
(here 
 $||\cdot ||$ stands for the distance to the nearest integer).
For irrational $\alpha$ put
\begin{equation}\label{la}
\lambda(\alpha) = \liminf_{t\to +\infty} t\cdot \psi_\alpha (t)
\end{equation}
and
\begin{equation}\label{de}
d(\alpha) = \limsup_{t\to +\infty} t\cdot \psi_\alpha (t).
\end{equation}

 {\it The  Lagrange spectrum} $\mathbb{L}$ is defined as
$$
\mathbb{L}=\{\lambda \in \mathbb{R}:\ \text{there exists} \
\alpha\in \mathbb{R} \ \text{such that} \
\lambda=\lambda (\alpha)\}.
$$
Here we should note that
\begin{equation}\label{gold}
\lambda\left(\frac{1+\sqrt{5}}{2}\right) = \frac{1}{\sqrt{5}}= 0.4472^+.
\end{equation}
is the maximal element of $\mathbb{L}$.
Also we note that Lagrange spectrum has a ``discrete part''
$$
 \frac{1}{\sqrt{5}},\frac{1}{\sqrt{8}},\ldots
$$
which is related to Markoff numbers and there exist a whole segment
$[0,\lambda^*]\subset \mathbb{L}$ which is known as {\it Hall's ray}.

 {\it The Dirichlet spectrum} $\mathbb{D}$  is defined as
$$
\mathbb{D}=\{d \in \mathbb{R}: \ \text{there exists} \
\alpha\in \mathbb{R} \ \text{ such that} \
d = d(\alpha )\}.
$$
 
The maximal element of $\mathbb{D}$ is $1$.
The minimal element from $\mathbb{D}$ is
$$
d\left( \frac{1+\sqrt{5}}{2}\right)
=\frac{1}{2} +\frac{1}{2\sqrt{5}}.
$$
It was found in \cite{Sze} (see also \cite{schmi}).
Dirichlet spectrum
 also has a ``discrete part''  and there exists a segment $[d^*, 1]\subset \mathbb{D}$.

 A lot of results related to Lagrange spectra one can find in \cite{cus}.
Dirichlet spectra was studied in \cite{divic}, \cite{japan}, \cite{i1}, \cite{i2}, \cite{i3},
An interesting survey one can find in \cite{malyshe}.

{\bf 2.Continued fractions.}
 
For irrational $\alpha$ represented as a continued fraction
 $$
\alpha=[a_0;a_1,a_2,\dots,a_t,\dots]=a_0+\frac{1}{{a_1+
\dfrac{1}{{a_2+\dots+\dfrac{1}{a_t+\dotsb}}}}}\,,\,\,\,\, a_0 \in \mathbb{Z},\,\, a_j \in \mathbb{N}.
$$
we define
$$
\alpha_\nu = [a_\nu;a_{\nu+1},...],\,\,\,\, \alpha_\nu^* = [0;a_\nu, a_{\nu-1},...,a_1].
 $$
We consider continued fraction's convergents
$$
\frac{p_\nu}{q_\nu} =[a_0;a_1,...,a_\nu].
$$
We should note that
\begin{equation}\label{beta}
\frac{q_{\nu-1}}{q_\nu} =[0;a_\nu,...,a_1].
\end{equation}

$$
\xi_\nu = ||q_\nu \alpha|| = |q_\nu \alpha - p_\nu|.
$$

It is a well-known fact that for irrational $\alpha$  the values of $\lambda(\alpha)$ and $d(\alpha)$ may be expressed as follows:
\begin{equation}\label{lala}
 \lambda(\alpha)
=\liminf_{\nu \to \infty} L(\alpha_{\nu}^*, \alpha_{\nu+1}),\,\,\,
\text{where}\,\,\,\, L(x,y) = \frac{1}{x+y},
\end{equation}
  \begin{equation}
\label{dede}
d(\alpha ) = \limsup_{\nu\to\infty}
D(\alpha_{\nu+1}^*,\alpha_{\nu+2}),\,\,\,\,\text{where}\,\,\,\,
D(x,y) = 
\frac{y}{x+y}.
\end{equation}

{\bf 3. Functions in two variables.}

In the sequel we consider two functions
$$
G(x,y) = \frac{x+y+1}{4}
$$
and
\begin{equation}\label{functor}
F(x,y) = \frac{(1 - xy)^2}{4(1+xy) (1-x)(1-y)} 
\end{equation}
in variables $x,y$.
It is clear that
\begin{equation}\label{clear}
 \min_{(x.y):\, x,y\ge 0} G(x,y)=\frac{1}{4},\,\,\,\,\,\
\max_{(x,y):\, x+y \le 1} G(x,y) = \frac{1}{2}.
\end{equation}

The following lemma is obvious.

  {\bf Lemma 1.}
\,\,{\it Put
$$
\Omega = \{(x,y)\in \mathbb{R}^2:\,\,\, 0\le x,  y < 1,\,\,  x+y^{-1} \ge 2, x^{-1}+y \ge 2\}.
$$
Then for the maximal value of $F(x,y) $ we have
$$
\max_{(x,y)\in \Omega}
F(x,y)  = \frac{1}{2}
$$
(the maximum here attains at any point satisfying
$x+y^{-1} =2$ or $ x^{-1}+y = 2$).
For the minimal value of $F(x,y)$ we have
$$
\min_{(x,y)\in \Omega}
F(x,y)  = F(0,0) =  \frac{1}{4}.
$$
}

{\bf 4. A function associated with Minkowski diagonal continued fraction.}

Given $\alpha \in \mathbb{R}$ we 
define a function 
$\mu_\alpha (t)$ which 
 corresponds to Minkowski diagonal continued fraction representation of $
\alpha$ (see \cite{MI}).  
To do this we recall  the Legendre theorem on continued fractions.
This theorem says that if 
\begin{equation}\label{1}
\left|
\alpha - \frac{A}{Q}\right|< \frac{1}{2Q^2},\,\,\,
(A,Q) =1
\end{equation}
then the fraction $\frac{A}{Q}$  is a convergent fraction for the continued fraction expansion of
$\alpha$. The converse statement is not true.
It may happen that
$\frac{A}{Q}$ is a convergent to $\alpha$
but (\ref{1}) is not valid.
We consider the sequence of the denominators of the convergents to $\alpha$
for which (\ref{1}) is true.
Let this sequence be
\begin{equation}\label{sece}
Q_0<Q_1<\cdots < Q_n<Q_{n+1}<\cdots .
\end{equation}
Then 
for $\alpha \not \in \mathbb{Q}$ the function
$\mu_\alpha (t) $ is defined by
\begin{equation}
\label{2}
\mu_\alpha (t) =
 \frac{Q_{n+1}-t}{Q_{n+1}-Q_n}\cdot ||Q_n \alpha||
+
 \frac{t-Q_n}{Q_{n+1}-Q_n}
\cdot ||Q_{n+1} \alpha||
,\,\,\
Q_n \le t \le Q_{n+1}.
\end{equation}
Here we should note that   for every $\nu$ one of the consecutive convergent fractions $\frac{p_\nu}{q_\nu},\,\frac{p_{\nu+1}}{q_{\nu+1}}$
to $\alpha$  satisfies (\ref{1}). 
So either
$$
(Q_n, Q_{n+1}) = (q_\nu, q_{\nu+1})
$$
for some $\nu$ and
$$
||Q_n \alpha || =\xi_\nu,\,\,\, ||Q_{n+1}\alpha|| =\xi_{\nu+1},
$$
 or
$$
(Q_n, Q_{n+1}) = (q_{\nu-1}, q_{\nu+1})
$$
for some $\nu$ and
$$
||Q_n \alpha || =\xi_{\nu-1},\,\,\, ||Q_{n+1}\alpha|| =\xi_{\nu+1},
$$
It is obvious that
$$
\liminf_{t\to +\infty} t\cdot \mu_\alpha (t) = \lambda (\alpha)
$$
whre the value $\lambda(\alpha)$ is defined in (\ref{la}).
Analogously to   $d(\alpha)$ defined in  (\ref{de}) we consider the value 
$$
\hbox{\got m} (\alpha) = \limsup_{t\to +\infty} t\cdot \mu_\alpha (t).
$$
We give here a result  analogous to formulas (\ref{lala}) and (\ref{dede}).

{\bf Theorem 1.}\,\,
{\it
Put
\begin{equation}\label{em}
\hbox{\got m}_n (\alpha)=
\begin{cases}
 G(\alpha_\nu^*,\alpha_{\nu+2}^{-1}),\,\,\,\text{if}\,\,\, (Q_n,Q_{n+1}) =(q_{\nu-1},q_{\nu+1})\,\,\,\text{with some } \,\,\,\nu,
\cr
F(\alpha_{\nu+1}^*,
\alpha_{\nu+2}^{-1}),\,\,\,\text{if}\,\,\, (Q_n,Q_{n+1}) =(q_{\nu},q_{\nu+1}) \,\,\,\text{with some } \,\,\,\nu.
\end{cases}
 \end{equation}
Then
$$
\hbox{\got m} (\alpha) = \limsup_{n\to +\infty} \hbox{\got m}_n (\alpha),
$$
}

We give a proof of Theorem 1 in Sections 6, 7, 8.

For example for
$$
\alpha = \frac{1+\sqrt{5}}{2}
=[1;1,1,1,...]
$$
we have $Q_n = q_n$ and
$$
\hbox{\got m} \left(\frac{1+\sqrt{5}}{2}\right) = F\left(\frac{\sqrt{5}-1}{2},\frac{\sqrt{5}-1}{2}\right) =\frac{1}{4}+\frac{1}{2\sqrt{5}}=0.4736^+,
$$
and for 
$$
\alpha = \sqrt{2}=[1;2,2,2,...]
$$
one has
\begin{equation}\label{qq}
\hbox{\got m}(\sqrt{2}) = F(\sqrt{2}-1,\sqrt{2}-1) =\frac{1}{4}+\frac{1}{4\sqrt{2}}=0.4267^+.
\end{equation}
But if we consider
$$
\alpha =\frac{1+\sqrt{3}}{2} = [1;2,1,2,1,2,1,...]
$$
then
\begin{equation}\label{ww}
\hbox{\got m} \left( \frac{1+\sqrt{3}}{2}\right) = G\left( \frac{\sqrt{3}-1}{2}, \frac{\sqrt{3}-1}{2}\right) = 
\frac{\sqrt{3}}{4}=0.4330^+.
\end{equation}
%Moreover if 
%$$
%\left|\alpha -\frac{p_\nu}{q_\nu}|\ge \frac{1}{2q_\nu^2}
%$$
%then $ a_{\nu+1} = 1

Alanogously to Lagrange and Dirichle spectra $\mathbb{L}$ and $\mathbb{D}$ we consider
the set
 $$
\mathbb{M} = \{ m \in \mathbb{R}: \,\,\, \exists \alpha \in \mathbb{R} \,\,
\text{such that}\,\, m = \hbox{\got m} (\alpha)\}.
  $$

{\bf Theorem 2.} \,\,{\it  For the minimal and the maximal element of $
\mathbb{M}$ one has
\begin{equation}\label{minimax}
\min \mathbb{M} = \frac{1}{4},\,\,\,\,\,\max \mathbb{M} = \frac{1}{2}.
\end{equation}
}

 We give a proof of Theorem 2 in Sections 9.

{\bf 5. Oscillating property.}

In \cite{KM}  it was proved that
for any two different irrational 
numbers
$\alpha , \beta$
such that $\alpha\pm \beta \not\in \mathbb{Z}$ the difference function
$$
\psi_\alpha (t) - \psi_\beta (t)
$$
changes  its sign infinitely many times as $t\to +\infty$.

The situation with oscillating property of the difference 
\begin{equation}\label{mudiff}
\mu_\alpha(t) - \mu_\beta(t)
\end{equation}
 is quite different.
In \cite{DARF} the 
it is shown that there exist real $\alpha $ and $\beta$
such that they are linearly independent over $\mathbb{Z}$ together with $1$ and
$$\mu_\alpha(t) > \mu_\beta(t),\,\,  \forall t \ge 1.$$
However  as it is shown in the same paper \cite{DARF} for {\it almost all} pairs $(\alpha,\beta)\in \mathbb{R}^2$
(in the sense of Lebesgue measure) 
the difference (\ref{mudiff}) does oscillate as $ t \to \infty$.

 The consideration of the values $\lambda(\alpha)$ and $\hbox{\got m}(\alpha)$ leads obviously to the following result.

{\bf Proposition 1.}\,\,{\it
Suppose that
real $\alpha, \beta$ satisfy
$$
\hbox{\got m}(\alpha ) < \lambda
 (\beta).
$$
Then there exists $t_0$ such that
$$
\mu_\alpha(t) < \mu_\beta (t),\,\,\,\, t\ge t_0.
$$}

In particular,
according to (\ref{gold},\ref{qq},\ref{ww})  one has
$$
\mu_{\sqrt{2}} (t) < \mu_{\frac{1+\sqrt{5}}{2}}(t),\,\,\,\,\,\,
\mu_{\frac{1+\sqrt{3}}{2}} (t) < \mu_{\frac{1+\sqrt{5}}{2}}(t)
$$
for all $t$ large enough.
This example shows that the conjecture  from \cite{DARF} is false.
Here we should note that J. Chaika \cite{CHA} was the first  to undestand that the conjecture 
from \cite{DARF} is false.

We would like to formulate here a positive result on oscillating of the difference (\ref{mudiff}).

{\bf Theorem 3.}\,\,{\it
Suppose that $\alpha$ and $\beta$ are quadratic irrationalities such that
they are linearly independent, together with $1$, over $\mathbb{Z}$ 
and
\begin{equation}\label{naturel}
 \lambda(\beta) < \lambda(\alpha) <\hbox{\got m}(\beta).
\end{equation}
Then the difference (\ref{mudiff}) changes its sign infinitely often as $t \to \infty$.}

We give a sketch of a proof for Theorem 3 in Section 10.

{\bf 6. Identities with continued fractions.}

{\bf Lemma 2.}
{\it

{\rm (i)}  
The following identities are valid:
\begin{equation}\label{cont1}
q_\nu \xi_{\nu} = \frac{1}{\alpha^*_\nu+\alpha_{\nu+1}} =  \frac{1}{(\alpha_{\nu+1}^*)^{-1}+(\alpha_{\nu+2})^{-1}}=
\frac{\alpha^*_{\nu+1}\alpha_{\nu+2}}{\alpha_{\nu+1}^*+\alpha_{\nu+2}},
\end{equation}
\begin{equation}\label{cont2}
\frac{\xi_\nu}{\xi_{\nu+1}} = \alpha_{\nu+2}.
\end{equation}

{\rm (ii)}
 Suppose that $a_{\nu+1} = 1$. Then
\begin{equation}\label{l1}
\frac{\xi_{\nu-1}}{\xi_{\nu+1}} = {\alpha_{\nu+2}+1}.
\end{equation}

}
Proof. The first equality from (\ref{cont1}) is well known (see \cite{schh}, Ch.1). To obtain the second one we should note that
$$
\alpha^*_\nu+\alpha_{\nu+2} =
\alpha^*_\nu+a_{\nu+1}
+
\frac{1}{\alpha_{\nu+2}} = \frac{1}{\alpha_{\nu+1}^*}+\frac{1}{\alpha_{\nu+2}}.
$$
Equality (\ref{cont1}) is proved.

To prove (\ref{cont2}) we observe that
$$
\frac{\xi_\nu}{\xi_{\nu+1}}=
\frac{q_\nu\xi_\nu}{q_{\nu+1}\xi_{\nu+1}} \, \frac{q_{\nu+1}}{q_\nu} =
(\alpha^*_{\nu+1}+\alpha_{\nu+2}) \times \frac{\alpha^*_{\nu+1}\alpha_{\nu+2}}{\alpha_{\nu+1}^*+\alpha_{\nu+2}}\times
\frac{1}{\alpha_{\nu+1}^*}
= \alpha_{\nu+2}
$$
(here we use (\ref{beta}) and equalities from  (\ref{cont1})).

To prove the statement from {\rm (ii)}
we observe that 
$$
\frac{\xi_{\nu-1}}{\xi_{\nu+1}} = \frac{\xi_{\nu-1}}{\xi_{\nu}}\, \frac{\xi_{\nu}}{\xi_{\nu+1}}
 =\alpha_{\nu+1}\alpha_{\nu+2}.
$$
Now we take into account  that
$$
\alpha_{\nu+1}=1+\frac{1}{\alpha_{\nu+2}},
$$
and (\ref{l1})
follows.$\Box$

Define 
$$
d_{1,\nu}= \sqrt{\frac{q_{\nu+1}-q_{\nu-1}}{\xi_{\nu-1}-\xi_{\nu+1}}},\,\,\,\,
d_{2,\nu}= \sqrt{\frac{q_{\nu+1}-q_{\nu}}{\xi_{\nu}-\xi_{\nu+1}}}
$$
and
$$
M_{1,\nu} =\frac{1}{4} \left( d_{1,\nu}\xi_{\nu-1}+\frac{q_{\nu-1}}{d_{1,\nu}}\right)^2,\,\,\,
M_{2,\nu} =\frac{1}{4} \left( d_{2,\nu}\xi_{\nu}+\frac{q_{\nu}}{d_{2,\nu}}\right)^2
.
$$

{\bf Lemma 3.}
\,\,{\it Suppose that $a_{\nu+1} = 1$. Then
\begin{equation}\label{c0}
M_{1,\nu}=
G(\alpha_\nu^*,\alpha_{\nu+2}^{-1}).
\end{equation}
}

Proof.
Calculations show that
\begin{equation}\label{c1}
d_{1,\nu}^2 \xi_{\nu-1}^2 =
\frac{q_\nu}{q_{\nu-1}} \frac{q_{\nu-1}\xi_{\nu-1}}{1-\frac{\xi_{\nu+1}}{\xi_{\nu-1}}}  
=\frac{(\alpha_{\nu+2}^{-1}+1)^2}{\alpha_\nu^*+ \alpha_{\nu+2}^{-1}+1}
\end{equation}
and
\begin{equation}\label{c2}
 \frac{q_{\nu-1}^2}{d_{1,\nu}^2} = 
\frac{q_{\nu-1}}{q_\nu} q_{\nu-1}\xi_{\nu-1} \left( 1- \frac{\xi_{\nu+1}}{\xi_{\nu-1}}\right) =
\frac{(\alpha_\nu^*)^2}{\alpha_\nu^*+ \alpha_{\nu+2}^{-1}+1}
\end{equation}
(here we use the equality (\ref{beta}), equality (\ref{cont1}) for $q_{\nu-1}\xi_{\nu-1}$ and (\ref{l1})).
Note that from   (\ref{cont1}) it follows that
\begin{equation}\label{c3}
q_{\nu-1}\xi_{\nu-1} =
\frac{1}{\frac{1}{1+\alpha_\nu^*}+\frac{1}{\alpha_{\nu+1}^{-1}}}
=\frac{\alpha_\nu^* (\alpha_{\nu+2}^{-1}+1)}{\alpha_\nu^*+ \alpha_{\nu+2}^{-1}+1}
.
\end{equation}
We combine (\ref{c1},\ref{c2},\ref{c3}) to get (\ref{c0}). Lemma is proved.$\Box$

%Proof. Put
%$$
%g = \frac{x-1}{1-y}
%.$$
%Then
 %$$
%\max_{(x,y)\in \Omega}
%F(x,y)
%=
%\frac{1}{4}\,
%\max_{g\ge 1} \left(\left( g+\frac{1}{g}+2\right)\times
%\max_{0\le y \le 1}\left( 
%\frac{1}{(1-g)y +g+1}\right)\right) = \frac{1}{2}.
%$$
%The last maximum attains at $g =1$.$\Box$ 

{\bf Lemma 4.}
\begin{equation}\label{c00}
M_{2,\nu}=
F(\alpha_{\nu+1}^*,
\alpha_{\nu+2}^{-1}).
\end{equation}

Proof. 
From the definition of $M_{2,\nu}$ and (\ref{cont1}) we have
$$
M_{2,\nu} =
\frac{1}{4((\alpha_{\nu+1}^*)^{-1}+(\alpha_{\nu+2})^{-1})}
 \,\frac{(g+1)^2}{g}
 ,\,\,\,\,
g =
\frac{\frac{q_{\nu+1}}{q_\nu}-1}{1-\frac{q_{\nu+1}\xi_{\nu+1}}{q_\nu\xi_\nu} \cdot \frac{q_\nu}{q_{\nu+1}}}.
$$
By
(\ref{beta}) and (\ref{cont2})
we get
$$
g =\frac{(\alpha_{\nu+1}^*)^{-1}-1}{1 -\frac{\alpha_{\nu+1}^*+\alpha_{\nu+2}}{(\alpha_{\nu+1}^*)^{-1}
+(\alpha_{\nu+2})^{-1}}\cdot\alpha_{\nu+1}^*}
=
\frac{(\alpha_{\nu+1}^*)^{-1}-1}{1-(\alpha_{\nu+2})^{-1}},
$$
and lemma follows by an easy calculation.$\Box$

{\bf 7. Inequalities with continued fractions.}

{\bf Lemma 5.}
\,\, {\it If $a_{\nu+1} = 1$
 then
\begin{equation}\label{ine2}
 d_{1,\nu}^2\cdot \frac{\xi_{\nu+1}}{q_{\nu+1}}
\le 1 \le d_{1,\nu}^2\cdot\frac{\xi_{\nu-1}}{q_{\nu-1}}.
\end{equation}
}

Proof. As $a_{\nu+1} = 1 $ by (\ref{beta}) and (\ref{l1}) we see that
$$
\frac{q_{\nu+1} -q_{\nu-1}}{q_{\nu+1}} = \frac{q_\nu}{q_{\nu+1}}=
\alpha_{\nu+1}^*\le 1<\alpha_{\nu+2}= \frac{\xi_{\nu-1}-\xi_{\nu+1}}{\xi_{\nu+1}}.
$$

The left inequality from  (\ref{ine2}) follows. To obtain the right inequality we should use 
the inequality 
$$ 1<
\frac{1}{\alpha_{\nu+1}^*}<
\frac{1}{\alpha_{\nu+1}^*}+\frac{1}{\alpha_{\nu+2}+1}
  $$
which  by (\ref{beta},\ref{cont2}) leads to
$$
\frac{q_{\nu+1} - q_{\nu-1}}{ q_{\nu-1}} =\frac{q_\nu}{q_{\nu-1}}
=
\frac{1}{\alpha_{\nu}^*}> 1-\frac{1}{\alpha_{\nu+2}+1} =
\frac{\xi_{\nu-1}-\xi_{\nu+1}}{\xi_{\nu-1
}}.
$$
This gives right inequality from (\ref{ine2}).$\Box$

{\bf Lemma 6. }
\,\, {\it Suppose that
$\alpha_{\nu}^*+\alpha_{\nu+1} >2$ and
$\alpha_{\nu+1}^*+\alpha_{\nu+2} >2$.
 Then
\begin{equation}\label{ine1}
 d_{2,\nu}^2\cdot \frac{\xi_{\nu+1}}{q_{\nu+1}}
\le 1 \le d_{2,\nu}^2\cdot\frac{\xi_{\nu}}{q_{\nu}}.
\end{equation}
}

Proof. 
From the conditions of lemma we have
$$
\frac{1}{\alpha_{\nu+1}^*}+\frac{1}{\alpha_{\nu+2}} =
\alpha_\nu^*+\alpha_{\nu+1} >2.
$$
So by (\ref{beta}) and (\ref{cont2}) we get
$$
\frac{q_{\nu+1}-q_\nu}{q_\nu} =
\frac{1}{\alpha_{\nu+1}^*}-1 > 1 -\frac{\xi_{\nu+1}}{\xi_\nu}
=\frac{\xi_\nu -\xi_{\nu+1}}{\xi_\nu}
$$
and hence
$$
\frac{\xi_\nu}{q_\nu} >
\frac{\xi_\nu - \xi_{\nu+1}}{q_{\nu+1}-q_\nu}.
$$
The last inequality coincides with the
right inequality from (\ref{ine1}).

To prove the left inequality from   (\ref{ine1}) we
observe that
$$
\frac{q_{\nu+1}-q_\nu}{q_{\nu+1}} = 1 -\alpha_{\nu+1}^* <
\alpha_{\nu+2} - 1 =
\frac{\xi_\nu}{\xi_{\nu+1}} - 1 =
\frac{\xi_\nu - \xi_{\nu+1}}{\xi_{\nu+1}}
$$
(here we use the condition 
$\alpha_{\nu+1}^*+\alpha_{\nu+2} >2$ and (\ref{beta},\ref{cont2}).$\Box$

 {\bf  8. Segments.}

In ${\mathbb R}^2$ with coordinates $(t,\mu)$ we consider  segments
$$
{\cal I}_\nu =[A_\nu, A_{\nu+1}],\,\,\,\,\,{\cal J}_\nu= [A_{\nu-1}, A_{\nu+1}],\,\,\,\,\,
A_j = (q_j,\xi_j)\in \mathbb{R}^2,\,\,\, j = \nu -1, \nu, \nu+1.
$$
 
{\bf Lemma 7.}\,\, {\it Suppose that  $a_{\nu+1} = 1$. Then
$$
\max_{(t,\mu) \in {\cal J}_\nu} t\cdot \mu = 
G(\alpha_\nu^*,\alpha_{\nu+2}^{-1}).
$$}

Proof. 
We consider 
 the hyperbolic rotation
$$
D_{1,\nu} = 
\left(
\begin{array}{cc}
 d_{1,\nu}^{-1}& 0\cr
0& d_{1,\nu}
\end{array}
\right).
$$
 This rotation  preserves each hyperbola $t \cdot \mu  = \omega$
and translates the segment ${\cal J}_\nu$ into the segment $[D_{1,\nu} A_{\nu-1},D_{1,\nu} A_{\nu+1}]$ which is orthogonal to the diagonal
$\{ \mu = t\}$
(this is the reason for the choice of the parameter $d_{1,\nu}$).
By Lemma 5
the endpoints of the segment $[D_{1,\nu} A_{\nu-1},D_{1,\nu} A_{\nu+1}]$
lie  by the different sides of the line $\{ t =\mu\}$.
So the
maximal value of the form $t\cdot \mu $ on the  segment occurs at the point
with $
\mu = t.
$
Easy calculation shows that  this maximum is equal to $M_{1,\nu}$. Now one should apply Lemma 3
and everything is proved.$\Box$
 
{\bf Lemma 8.}
\,\, {\it Suppose that
$\alpha_{\nu}^*+\alpha_{\nu+1} >2$ and
$\alpha_{\nu+1}^*+\alpha_{\nu+2} >2$.
Then
$$
\max_{(t,\mu) \in {\cal J}_\nu} t\cdot \mu = 
F(\alpha_{\nu+1}^*,
\alpha_{\nu+2}^{-1})..
$$}

Proof.  The proof is similar to the proof of Lemma 7. One should consider
 the hyperbolic rotation
$$
D_{2,\nu} = 
\left(
\begin{array}{cc}
 d_{2,\nu}^{-1}& 0\cr
0& d_{2,\nu}
\end{array}
\right)
$$
applied to the segment 
${\cal  I}_\nu$ and take into account Lemmas 4, 6.$\Box$

Now we should note that Theorem 1 immediately follows from the definition of $\hbox{\got m}(\alpha)$
and Lemmas 7,8.

{\bf 9. Proof of Theorem 2.} 

To prove (\ref{minimax}) we use
equalities (\ref{clear}) for $G(x,y)$, Lemma 1  and Theorem 1.
This gives
$$
\min \mathbb{M} \ge \frac{1}{4},\,\,\,\,\,\max \mathbb{M} \le \frac{1}{2}.
$$
To prove that here are just equalities one should consider examples
$$
\alpha^- = [0;a_1,a_2,...,a_n,...]
,\,\,\,\,\, a_n \to \infty,\,\, n \to \infty
$$
and
$$
\alpha^+ = [0;1,1,a_3,1,1,a_6,1,1,a_9,1,1,a_{12},1,1,a_{15},...],\,\,\,\,\,
a_3<a_6<a_9<a_{12}<a_{15}<...
$$
with
$$
\hbox{\got m}(\alpha^-)
=\lim_{\nu\to \infty}F(\alpha_{\nu+1}^*,
\alpha_{\nu+2}^{-1})  = \frac{1}{4}
$$
and
$$
\hbox{\got m}(\alpha^+)
=\lim_{\nu\to \infty}F(\alpha_{3\nu+2}^*,
\alpha_{3(\nu+1)}^{-1})
=
\lim_{\nu\to \infty}G(\alpha_{3\nu+1}^*,
\alpha_{3(\nu+1)}^{-1})
  = \frac{1}{2}.
$$

{\bf 10. Proof of Theorem 3.}

From the conditions of Theorem 3 we see that $
 \beta \not\in \mathbb{Q}(\alpha)$.
Consider fundamental units $\varepsilon_\alpha$ and  $\varepsilon_\beta$
of the fields $\mathbb{Q}(\alpha)$ and $\mathbb{Q}(\alpha)$,
respectively. Then
\begin{equation}\label{irra}
\frac{\log \varepsilon_\alpha}{\log \varepsilon_\beta} \not \in \mathbb{Q}.
 \end{equation}
Given positive $\eta$ one can construct
three sequences of integers
$$
Q_\nu^{[\beta],1},\,\,\,\, 
Q_\nu^{[\beta],2},\,\,\,\, 
 Q_\nu^{[\alpha]},\,\,\,\, \nu = 0,1,2,3,...
$$
such that
$$
|
 Q_\nu^{[\beta],1} \mu_\beta (Q_\nu^{[\beta],1}) - \lambda({\beta})|
<\eta
$$
and 
$$
 Q_\nu^{[\beta],1} 
= C_{\beta,1} B_1^\nu (1+o(1)), \,\,\, \nu \to +\infty,
$$
where $C_{\beta, 1}$ is positive and $B_1 = \varepsilon_\beta^{b_1} $ with a positive integer $b_1$;
$$
|
 Q_\nu^{[\beta],2} \mu_\beta (Q_\nu^{[\beta],2}) - \hbox{\got m}({\beta})|
<\eta
$$
and 
$$
 Q_\nu^{[\beta],2} 
= C_{\beta,2} B_2^\nu (1+o(1)), \,\,\, \nu \to +\infty,
$$
where $C_{\beta, 2}$ is positive and $B_2 = \varepsilon_\beta^{b_2} $ with a positive integer $b_2$;
$$
 | Q_\nu^{[\alpha]} \mu_\alpha (Q_\nu^{[\alpha]}) - \lambda{(\alpha)}|
<\eta
$$
and 
$$
  Q_\nu^{[\alpha]} 
= C_{\alpha} A^\nu (1+o(1)), \,\,\, \nu \to +\infty,
$$
where $C_{\alpha}$ is positive and $A = \varepsilon_\alpha^{a} $ with a positive integer $a$.

From (\ref{irra}) we have
$$
\frac{\log{B_j}}{\log A} \not\in \mathbb{Q},\,\,\, j = 1,2.
$$
By the Kronecker theorem given real $\omega_1, \omega_1$ and positive $\eta$ there exist infinitely many
pairs $(\nu_{1,n}, \kappa_{1,n})$ and  $(\nu_{2,\nu}, \kappa_{2,n})$ such that
$$
 |\nu_{1,n}\log B_1 - \kappa_{1,n}\log A - \omega_1 | <\eta,
$$
and
$$
 |\nu_{2,n}\log B_2  - \kappa_{2,n}\log A - \omega_2 | <\eta,
$$
With the proper choice of $\omega_1,\omega_2$ this gives
$$
\frac{Q_{\nu_{1,n}}^{[\beta],1} }{  Q_{\kappa_{1,n}}^{[\alpha]}} \to 1,\,\,\,\,
\frac{Q_{\nu_{2,n}}^{[\beta],2} }{  Q_{\kappa_{2,n}}^{[\alpha]}} \to 1,\,\,\,\,
n\to \infty.
$$
Simple calculation shows that
for small $\eta$ and for $n$ large enough one has
$$
\left|
\frac{
  \mu_\beta (Q_{\nu_{j,n}}^{[\beta],j})  
}
{
  \mu_\beta (Q_{\kappa_{j,n}}^{[\alpha]})
 }  - 1
\right| < 2\eta,\,\,\,\, j = 1,2.
$$
and

So for $n$ large enough
$$
\frac{\mu_\beta (Q_{\kappa_{1,n}}^{[\alpha]})}
{\mu_\alpha (Q_{\kappa_{1,n}}^{[\alpha]})}<
\frac{\lambda (\beta)}
{\lambda (\alpha)}\, (1+4\eta)<1<
\frac{\hbox{\got m}(\beta)}
{\lambda (\alpha)}\, (1-4\eta)
<
\frac{\mu_\beta (Q_{\kappa_{2,n}}^{[\alpha]})}
{\mu_\alpha (Q_{\kappa_{2,n}}^{[\alpha]})}
$$
The last inequalities show that the difference (\ref{mudiff}) does oscillate.$\Box$

\end{document}